\documentclass[12pt,]{article}
\usepackage{amsmath,amsthm}
\usepackage{amssymb}

\usepackage{graphicx}

\usepackage{color} 

\usepackage{setspace}
\makeatletter
\renewcommand{\@biblabel}[1]{\quad#1.}
\makeatother

\pdfoutput=1

 \newtheorem{thm}{Result}

\newtheorem{rem}{Remark}

\def\ind{{\mathchoice {\rm 1\mskip-4mu l} {\rm 1\mskip-4mu l}
{\rm 1\mskip-4.5mu l} {\rm 1\mskip-5mu l}}}

\begin{document}

\begin{flushleft}
{\Large
\textbf{Mean-Field Games for Marriage}
}
\\
Dario Bauso$^{1}$, 
Ben Mansour Dia$^{2}$,
Boualem Djehiche$^{3}$,
Hamidou Tembine$^{2,\ast}$, 
Raul Tempone$^{2}$
\\
\bf{1} Dipartimento di Ingegneria Chimica, Gestionale, Informatica, Meccanica, Palermo, Italy
\\
\bf{2} SRI Center for Uncertainty Quantification in Computational Science and Engineering, 
							King Abdullah University of Science and Technology
\\
\bf{3} Department of Mathematics, Royal Institute of Technology (KTH), Stockholm, Sweden
\\
$\ast$ E-mail: tembine@ieee.org
\end{flushleft}

\section*{Abstract}

This article examines mean-field games for marriage.  The  results support the argument that optimizing the long-term well-being  through effort and  social  feeling state distribution (mean-field)  will help to  stabilize marriage.  However ,  if the cost of effort is very high, the couple fluctuates in a bad feeling state or the marriage  breaks down.  We then examine the influence of  society on a couple using mean field sentimental games.  We show  that, in mean-field equilibrium,  the optimal effort is always higher than the one-shot optimal effort.  We illustrate numerically  the influence of the couple's network on their feeling states and their well-being.

\section*{Significance}
The myth of marriage has been and is still a fascinating historical societal phenomenon. Paradoxically, the empirical divorce rates are at an all-time high.  In order to design and evaluate interventions, theoretical understanding of marital stability and dissolution is crucial. 
This article describes a unique paradigm for preserving relationships and marital stability using mean-field game theory. 
The approach is quite general for modeling social interaction, and can be applied to empirical data generated over time. 
\section*{Introduction}

{\color{black}
We look at marital interactions and relationships. Here relationship refers to the unit (couple), rather than to the two persons. We focus on a very ephemeral thing:  what happens to the couple as they interact over time.
It is not either person, it is something that happens when they are together, i.e., in a couple. The couple will create something that we call a "feeling state"  as they talk to each other, as they smile, as they move. We model the feeling state of a couple as a noisy differential equation influenced by random events, parents, children, friends and the social distribution of marital status (mean-field of states), where the term "mean-field" refers to the distribution of couple feeling states in the society. 
}

\noindent {\bf Why a mathematical model for marital interaction?}
Putting the study of social relationships on  mathematical footing represents a major advance in our ability to understand and perhaps  regulate these relationships for the betterment of all mankind \cite{refty0}. The novelty of the present study inheres in the analysis  of the social distribution of feeling states (mean-field) on a generic couple.

\noindent{\bf On the influence of couple-feeling-states in society :} Since the society, friends and parents  of a couple may influence them in an aggregative manner, a mean-field approach  is suitable for such an interaction. Here the distribution of states (mean-field)- of all couples in the society plays an important  role in the payoff and in the effort used to maintain a marriage. Game theory is a branch of mathematics which studies strategic interactions. Our idea is to use  mean-field games to capture the influence of society distribution of states on a generic couple.

In their pioneering work in \cite{refty0}, Gottman et al.  have widely  illustrated the importance of mathematical theory  in the field of marital research. In particular, they have shown that mathematical models help us to better understand the sentimental dynamics and hence to propose an appropriate intervention or  therapy for a couple. The authors have developed a model to analyze marital instability, by calibrating a system of difference equations for the evolution of partners' emotions during a conversation. The interaction dynamics is described by a function of experience and whose intuitive
understanding of the influence among partners during conflict. More recently, Rey \cite{senti} formulated  the sentimental relationship of a couple as an optimal control problem. A state variable monitors the wellness of the relationship whose natural decay in time must be counteracted with effort  according to a widely accepted principle in marital psychology. 
Stationary sentimental steady states are examined in \cite{senti,senti2}  the context of infinite horizon
 discounted optimal control problems.

 None of these previous works considered the influence of model uncertainty, noise, random events, and the 
social distribution of states (mean field) . In this work we examine these  important issues.

This paper has three main features: 
\begin{itemize}
\item[(i)] We introduce stochastic optimal control in marital interaction. We explicitly characterize the optimality equations under random fluctuations which comes from model uncertainty in the drift function of the feeling state.
\item[(ii)] We examine the impact of mean-field and effort in the  feeling state of the couple. 
\item[(iii)] We propose a mean-field game model for marriage in a non-linear setup where generically two stable 
states (divorce and marriage) are observed. We study the impact of the network of couples on the feeling state. We show that, in  mean-field equilibrium,  the optimal effort of the mean-field sentimental game is always higher than the one-shot optimal effort if the social contribution is controlled.
\end{itemize}
While the study here focuses on marital interaction, the mathematical methodology can be applied to other situations such as  a society's ecological or smoking behavior. To the best our knowledge this is the first attempt  to model and understand marital interactions through mean-field influence and stochastic games.

The remainder of the paper is organized as follows. In the next section we present a generic model of  sentimental dynamics examining  both controlled and uncontrolled situations. After that we introduce mean-field sentimental games. Then, we analyze the optimal behavior and mean-field sentimental equilibria with and without noise. As an illustration we examine an interesting sentimental dynamics wherein marriage and divorce can be seen as ''stable" steady states. Adding  strategic effort behavior into such couple  dynamics  can change the long-term behavior of their feeling state if the system is uncontrolled. 

%
%
%
%
\section*{Model and Methods}
{\color{black} We first present a basic mathematical framework for feeling state and well-being of a couple. We will explain below how the model can be modified to include the (positive or negative) impact of the society (and the "social pressure") on a generic couple.  }
We describe the sentimental dynamics, define the payoff functional for the marital interaction and state the optimal control problem. 
\subsection*{Governing feeling state equation }
The feeling state (feeling level) $x$ of the couple  is modeled as an It\^o's stochastic differential equation:
\begin{eqnarray}
\label{modelrey} 
\displaystyle{ dx = f(x,e)dt + \sigma d{B}(t),\;\;\text{for}\;\; t\in [0,T],}
\end{eqnarray}
where $f(x,e) = - h(x)+ a e$,  $ a>0$ represents the expected variation of the feeling in a short time window, $\sigma \geq 0$, $h$ is  smooth and goes to infinity with $x.$  ${B}(t)$ is a standard Brownian motion. The effort function $e(t)$ can be controlled by the couple. The sentimental dynamics is  expected to start at a high feeling level $x_0$ greater than $\underline{x}$ where $\underline{x}$ is a certain threshold value below which the relationship of the couple is not considered as satisfactory. The parameter $a$ represents the efficiency of the effort. 
\begin{rem} The deterministic part of the governing feeling state is an experimental derivation from marital studies \cite{refty0,senti,senti2}. {\color{black} However, here, the function $h$ is not necessarily monotone in $x.$ }
We introduce a stochastic term into the sentimental dynamics of the existing literature \cite{refty0,senti,senti2} for multiple reasons: {\color{black} (1)  shocks, (2) random events in the social network of the couple may affect the couple, (3) the drift funtion $f$ may not be perfectly known. We add an uncertainty term into it. }
\end{rem}
\begin{rem}
 If the marriage  starts with a feeling level $x_0$ and is stopped at the first time that $x(t)$ is below a certain threshold $ \underline{x}-\epsilon$  or one of the couple member dies (life expectancy is set to $T_0 = 100$ years)  then the length of the horizon is expected to be finite. This is why we consider a finite horizon problem. Let $T_{\mbox{divorce}}=\inf\{t >0 |  \ x(t) \leq \underline{x}-\varepsilon \}$, one can take $T = \min (T_0, T_{\mbox{divorce}})$ which means that $T$ is a random variable. 
\end{rem}
\subsection*{Setting of the problem}
\noindent Based on the feeling state and the effort, we define the payoff functional for a couple during 
$\ [0,T], \ T>0$ as
\begin{eqnarray}
  P(x_0 , e) = \mathbb{E}\left(g(x(T)) + \int_0^T [s(x(t))- c(e(t))] dt \right),
\end{eqnarray}
where $s(x(t))$ represents the well-being of the couple at state $x(t)$ and $c(e(t))$ is the cost associated to the effort $e(t).$
The instantaneous payoff captures the risk at time $t$ but we have limited our analysis to the long-term risk-neutral setup. More details on risk-sensitive mean-field type control and games can be found in \cite{rsmfg,boualemtembine}. 

\noindent We define a couple as a union of  two adult persons having certain emotional and physical interactions, 
and living  together. We define family as the direct  ascendants and the direct descendants of the conjoint. In this sense, we read the quantity $g(x(T))$ as the heritage that the couple bequeaths to the family in terms of well-being. This heritage  will influence the next generation sentimental game.  It is reasonable to assume that the function $g$ increases with the 
feeling state and $g'(x)\geq 0$.

\noindent The function $s$ is the well-being function of the couple. It depends on the feeling state $x$ in the
 sense that, the higher  the feeling state, the more joyful the couple. The implicit purpose of the control problem is 
to increase well-being while taking into consideration the cost. Thispurpose achieved by expanding the feeling state above a certain satisfactory level. 

\noindent Contrastingly, couples with low feeling state are more or less happy and have the tendency to seize again. This tendency
 is expressed by taking actions; caring for children, giving gifts, cooking special meals, making concessions, in short, handling the couple according to his or her expectation. We define the effort as a constructive action from one conjoint to the other conjoint. It might be useful for the reader to notice that our definition of effort does not take into  account  counterproductive actions. For example, knowing that Juliette does not like football, it is negative that Romeo takes her to the Camp Nou Stadium  for 
FC Barcelona vs Milan AC making her miss her preferred telenovella movies. The cost of effort  can represent  the psychological, financial and/or emotional load provided by the couple when accomplishing this effort. Clearly, greater effort is more expensive as it imposes a higher cost. It is appropriate then to say that the couple  continuously provides efforts to achieve a greater well-being of the couple itself.

\begin{figure}[htb]
\centering
\includegraphics[width=0.4\textwidth]{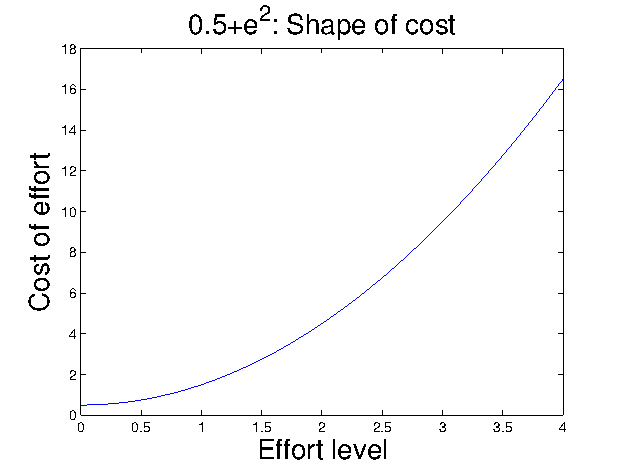}
\caption{Typical shape of the cost function.}
\label{fig1:plotwellbeing}
\end{figure}
\begin{figure}
\centering
\includegraphics[width=0.4\textwidth]{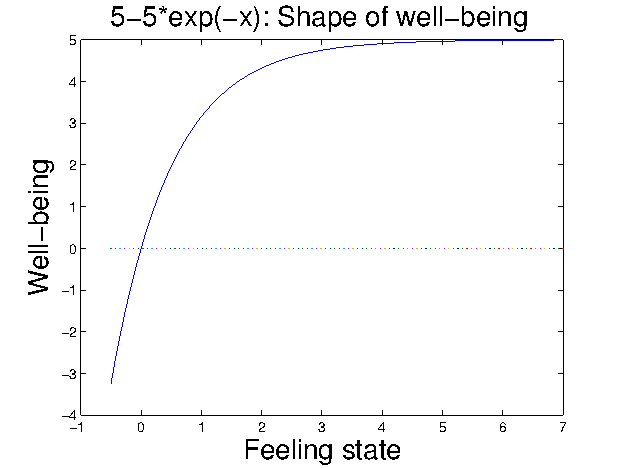}
\caption{Typical shape of  well being.}
\label{fig2:plotwellbeing}
\end{figure}

We adopt the following assumptions:  
\begin{itemize}
\item the drift function $\ f: \ \mathbb{R}\rightarrow \mathbb{R}$ is continuously differentiable ($\mathcal{C}^1$), 
\item the  well-being function  $s: \ \mathbb{R}\rightarrow \mathbb{R}$ is $\mathcal{C}^1$-differentiable, non-decreasing, concave and saturated at $+\infty.$,
\item the cost function $c: \ \mathbb{R}\rightarrow \mathbb{R}$  is twice continuously differentiable ($\mathcal{C}^2$), non-decreasing and strictly convex.
\end{itemize}

In  Figure  \ref{fig1:plotwellbeing}-\ref{fig2:plotwellbeing}, we illustrate the typical well-being and cost functions.

\noindent  In the next sections, we consider respectively an open-loop control problem, an optimal feedback strategy for the couple and a mean-field sentimental strategy.
%
%
%
%
%
%
\section*{Results }
\subsection*{Optimal open-loop effort control}
A control law of the form $e(t) = e(t; t_0, x_0)$ for $ t\in [t_0, T]$, determined for a particular initial state value $x(t_0) = x_0$ is called an open-loop effort. This does not depend explicitly on the state $x(t)$. One seeks for an optimal open loop effort strategy by applying the maximum principle to the following problem
\begin{eqnarray}
\label{open_loop}
(\mathcal{P}_0) 
\left\{
\begin{array}{lll}
\displaystyle{ \sup_{(e(t))_t} P(x_0, e)}, \\
\displaystyle{\mbox{ subject to } (\ref{modelrey})\ \mbox{ that starts at}\  x_0 > \underline{x} + \varepsilon.}
\end{array}
\right.
\end{eqnarray}
To do so, we use the Hamiltonian and the adjoint process. First we analyze the deterministic case ($\sigma = 0$) and  then the open-loop noisy case.
%
%
\subsubsection*{Deterministic open-loop optimal effort}
When $\sigma=0$, the Hamiltonian is  
\begin{eqnarray}
H_0(x,p,e) = f(x,e)  p + s(x)-c(e).
\end{eqnarray}
A maximizer of $H_0(x,p,e)$ with the respect to the effort provides an open-loop optimal control associated with the co-state (adjoint) process $p$.

\noindent The adjoint process of the optimal control is
\begin{eqnarray} 
\label{costate_open_loop}
 p_t = -H_{0,x} dt = (- p f_x -s'(x))dt, 
\end{eqnarray}
with $p(T) = g'(x_T)$. Note that $H_0$ is strictly concave  in $e$. For concave function, the first order optimality  condition is also a sufficient condition for interior point. The (interior) effort should 
satisfies $c'(e)=ap\geq 0.$
The positivity constraint of the effort suggests $p\geq 0$ and then there is a unique solution.

\noindent The open-loop optimal control system via Pontryagin maximum principle \cite{yong} yields:
\begin{eqnarray}
\label{opt_open_loop}
\left\{
\begin{array}{lll}
\displaystyle{e^*(t)   = \max(0, (c')^{-1}[a p(t)]),}
\vspace{0.1cm}\\
\displaystyle{\dot{p}(t) = (p(t) h'(x(t)) - s'(x(t))),}
\vspace{0.1cm}\\
\displaystyle{\dot{x}(t) = -h(x(t))+ \max(0, a (c')^{-1}[a p(t)]),}
\vspace{0.1cm}\\
x(0)=x_0,\ \ p(T)=g'(x_T).
\end{array}
\right.
\end{eqnarray}
\noindent By differentiating $e^*(t)$ and combining with the above system, we arrive at
\begin{eqnarray}
\dot{e}  = \displaystyle{\frac{1}{c''(e)}\left[ c'(e)h'(x)-a s'(x)\right]\ind_{\{ e(t)\geq 0\}}.}
\end{eqnarray}
Hence, one gets the following dynamical system between the optimal control and the optimal feeling:
\begin{eqnarray}
\label{opt_open_loop_dyn}
\left\{
\begin{array}{lll}
\displaystyle{x(0)    =  x_0,}\\
\displaystyle{\dot{x} = -h(x)+a e,}\\
\displaystyle{\dot{e} = \frac{1}{c''(e)}\left[ c'(e)h'(x)-a s'(x)\right]\ind_{\{ e(t)\geq 0\}},}\\
\displaystyle{e(T)    = \max\left(0, (c')^{-1}[a g'(x(T))]\right).}
\end{array}
\right.
\end{eqnarray}

Our first question is about the  well-posedness of the above system (\ref{opt_open_loop_dyn}). Our goal is to provide a sufficiency condition for existence of a solution. 
We use fixed-point theorem to establish existence under the above assumptions. Since the functions $s$ and $h$ are $\mathcal{C}^1$-differentiable and the  strictly convex function $c$ is $\mathcal{C}^2$-differentiable, (\ref{opt_open_loop_dyn}) admits a local solution.  Next we  use a stochastic maximum principle approach to  analyze the case where $\sigma$ is nonzero, i.e., the stochastic case.  
%
%
%
%
\subsubsection*{Stochastic optimal open-loop effort}

Following the stochastic maximum principle technique, the  adjoint processes for constant variance coefficient  yields
\begin{eqnarray} 
\label{sto_open_loop_adjoint}
\left.
\begin{array}{lll}
\displaystyle{dp(t)= (-p(t)  f_x - s'(x(t)))dt + q  d{B}(t)} 
\vspace{0.1cm}\\ 
\text{with} \;\; p(T) =  g'(x(T)),
\end{array}
\right.
\end{eqnarray}
where $q$ is the adjoint variable associated with the diffusion term. Let $\tilde{c}$ be the Legendre-Fenchel transform of $c:$
\begin{eqnarray}
\tilde{c}(ap(t))=\sup_{e}\left\lbrace ap(t) \cdot e- c(e)\right\rbrace. 
\end{eqnarray}
It is well known that $\tilde{c}(ap(t))$ is convex in $p(t)$. Furthermore, the optimal control is given by 
\begin{eqnarray}
e^*(t)=\max(0, \tilde{c}'(ap(t))).
\end{eqnarray}
Using Ito's calculus, the Euler-Lagrange system is given by
\begin{eqnarray}
\label{stoc_opt_open_loop_dyn}
\left\{
\begin{array}{lll}
\displaystyle{x(0)    =  x_0,}
\vspace{0.1cm}\\
\displaystyle{dx(t)   = (-h(x(t))+a e(t))dt+\sigma d{B}(t),}
\vspace{0.1cm}\\
\displaystyle{de(t)   = \left[ \tilde{c} ''(c'(e(t)))  \left(c'(e(t)) h'(x(t))
                       - as'(x(t))\right) \right. } 
											\vspace{0.1cm}\\   
\displaystyle{  \hspace{1cm}  \left. + \frac{q^2}{2}a^2 \tilde{c} '''(c(e(t))) \right] dt + \tilde{c} ''(c'(e(t))) aqd{B}(t),}
\vspace{0.1cm}\\ 
\displaystyle{e(T)    = \max\left(0, (c')^{-1}[a g'(x(T))]\right).}
\end{array}
\right.
\end{eqnarray}

A sufficiency condition for existence of solution $(x(t),e(t))$ is
 the uniformly Lipschitz with respect to $(x,e)$ of the coefficient functions and uniformly Lipschitz condition
 for $(c')^{-1}$ over the horizon $[0,T].$

\subsection*{Optimal feedback effort strategy  for the couple}
A feedback effort law is of the form $e(t, x(t), t_0, x_0), \; t\in [t_0, T],$ i.e., the effort depends on time $t$, feeling $x(t)$ and possibly the initial conditions. Our motivation for feedback strategy comes from the following result from Jackson (1957), page 79: ''\textit{ A family interaction is a closed information system in which the variation in outputs or behavior are feedback in order to correct the system's response}''. This idea was re-used in the book \cite{refty0} entitled ''The  Mathematics of Marriage".
The pay-off functional is  written as
\begin{eqnarray}
  P_{\text{\tiny{closed}}}(x,e)=\mathbb{E}\left(g(x(T))+\int_t^T \left[s(x({t'}))-c(e({t'}))\right] dt' \ | \ x(t) = x\right).
\end{eqnarray}
To find the optimal feedback control $e^*(t,x),$ we use dynamic programming. Let $v(t,x)$ be the value starting from $x$ at time $t,$ i.e.,
\begin{eqnarray}
\label{opt_feedback_problem}
(\mathcal{P}_1)
\left\{
\begin{array}{lll}
\displaystyle{v(t,x)=\sup_{e(t,.)\geq 0} P(x,e),}
\vspace{0.1cm}\\
\displaystyle{\text{subject to the feeling state dynamics (\ref{modelrey}).}}
\end{array}
\right.
\end{eqnarray}
The Hamiltonian is given by
\begin{eqnarray}
\label{opt_feedback_hamiltonian}
H(x,p) = - h(x)p + s(x) + \tilde{c}(a p). 
\end{eqnarray} 
{\color{black} Since $\tilde{c}(a p)$ is convex $p$ we deduce that  $H$ is convex in $p$.} It is well known that, if there exists a twice continuous differentiable function $v(t,x)$  solution of 
\begin{eqnarray}
\label{opt_feedback_cost}
\left\{
\begin{array}{lll}
 v_t -h(x) v_x + s(x) + \tilde{c}(av_x) + \frac{\sigma^2}{2}v_{xx} = 0 \\
 \text{with}\;\; v(T,x)= g(x),
 \end{array}
\right.
\end{eqnarray}
then the stochastic optimal control effort is
\begin{eqnarray}
\label{opt_feedback_effort}
\displaystyle{e^*(t,x) = \max\left(0\, , \, \tilde{c}'(av_x)\right).}
\end{eqnarray} 
 Under the assumptions on $c, h, s, g$ and given $x_0$ and the strictly concavity of the payoff function in $e$ implies that there is a unique optimal control. Thus, the existence of the value  follows.
%
%
%
%
\subsection*{Mean field sentimental games } \label{mfg}
There is a significant consensus around the idea  that the society  state (particularly the feeling states of the friends, parents, friends' of friends and social environment) may influence the status of the couple. For example, if many of the friends are in lower feeling states, their divorce can be contagious. Moreover, a couple's tendency to divorce depends not just on their friends' divorce status, but also extends to their friend's friends and so on. Thus we must  consider the distribution of states within the entire society, i.e., the mean-field of states. Previous works neglected the influence of society and the couple's network \cite{senti,senti2} . In this paper, we take into consideration the mean-field state of the society in the sentimental dynamics of the couple. This leads to a {\it mean-field sentimental game}. 

\noindent Roughly speaking, we introduce a society's social influences (mean-field) and  external shocks into the marital interaction model. To do this, we assume that the function $h$ has the form $\hat{h}(x(t),m(t))$ where $m(t,.)$ is the distribution of states of all couples in the society at time $t$. 
Based on the model (\ref{modelrey}), we state that the feeling state (feeling level) $x(t)$ of the couple  is modeled as a McKean-Vlasov It\^o's stochastic differential equation:
\begin{eqnarray}
\label{modelbook2} 
\displaystyle{ dx(t) = \left[- \hat{h}(x(t),m^*(t)) + a e(t)\right] dt  + \sigma d{B}(t).}
\end{eqnarray} 
{\color{black} where $m^*$ is the equilibrium distribution of states.}

In this sequel the long-term payoff of a generic couple is written as 
\begin{eqnarray}
\label{mf_payoff}
\displaystyle{  P_{\text{\tiny{field}}}(x_0, m^*, e)= \mathbb{E}\left(g(x(T)){+}\int_0^T \left[\hat{s}(x(t), m^*(t)){-}c(e(t))\right] dt \right).}
\end{eqnarray}

The objective of each generic couple is to maximize its long-term payoff through the mean-field fixed-point problem:
\begin{eqnarray}
\label{mean_field}
(\mathcal{P}_2)
\left\{
\begin{array}{lll}
\displaystyle{ \sup_{(e(t))_t} P_{\text{\tiny{field}}}(x_0, m^*, e)}, \\
\displaystyle{\mbox{subject to the feeling state dynamics }(\ref{modelbook2}).}
\end{array}
\right.
\end{eqnarray}
The happiness function $\hat{s}$ covers the well-being of the couple and the satisfaction of its network. It  is natural that the happiness distribution over the network is concentrated in the family set. It is also known that the social network  may contribute to enhance the well-being of the couple.

\subsubsection*{Mean-field sentimental equilibria }
We define a mean-field sentimental equilibrium as a (Nash/Wardrop) equilibrium of the mean-field sentimental game. A mean-field sentimental equilibrium in feedback strategies is a situation in which no couple has incentive to move unilaterally from its effort feedback strategy.

\noindent Following \cite{mfg2007}, the mean-field equilibria are solution of the following backward-forward system  for $(t,x) \in  [0,T]\times \mathbb{R}$ is
\begin{eqnarray}
\label{mf_system}
\left\{
\begin{array}{lll}
\displaystyle{  v_t - \hat{h}(x,m) v_x  + \hat{s}(x,m)  +  \tilde{c}(a  v_x)  + \frac{\sigma^2}{2}v_{xx} = 0  ,}\\
\displaystyle{ m_t + \partial_x \left[- m \hat{h}  + a m \tilde{c}'\left(a  v_x)\right)\right]  - \frac{\sigma^2}{2}m_{xx}= 0,}\\
\displaystyle{m(t = 0 , x) =  m_0(x)\;\; \text{for }  \ x \in \mathbb{R},} \\ 
 \displaystyle{ v(t = T , x) = v_T(x) = g(x)\;\; \text{for } \ x \in  \mathbb{R}.}\\
\end{array}
\right.
\end{eqnarray}

Under the above assumptions on $\hat{h},\; \hat{s},\; c$ and the set of effort strategies, the mean-field sentimental game has an equilibrium.

\begin{thm} Asssume that the contribution of the mean-field to the feeling state of a generic couple is small. Then,
the optimal level of effort for a long-term viability of the couple (during their lifetime)  which keeps a happy relationship going is always greater than the effort level that would be chosen in a one-shot, i.e., 
$e^*(t,x)\ \geq \underline{e} $, where $c'(\underline{e})=0$.
\end{thm}

By assumption, the functions $c'$ and $g$ are non-decreasing, and the term $h_2(m)$ has only small influence in the drift. By small influence of $h_2(m)$ we mean that the uncontrolled system leads to a divorce or fluctuate feeling state if no effort is injected.
The optimality equation (\ref{mf_system}) and the strict convexity of the cost $c$ gives $c'(\underline{e})=0\leq a v_x(t,x)=c'(e^*(t,x)).$ This means that 
$c'(\underline{e})\leq c'(e^*(t,x)).$ Using the inverse function $(c')^{-1},$ one gets $\ \underline{e}\leq e^*(t,x)$ and there are several time $t\in [0,T]$ for which $\ \underline{e}< e^*(t,x).$

%
%
\subsubsection*{Open-loop sentimental equilibria}
\noindent In this section, we analyze the influence of the society state (mean field) on the couple.  Besides, we assume that the functional $\hat{h}$, $\hat{s}$ and $\hat{c}$ take on the following forms:
\begin{eqnarray*} \label{mf_h_and_s}
\hat{h}(x,m) = h(x) - h_2({m}),\ \hat{s}(x,m) = s(x) + s_2({m}), 
\end{eqnarray*}
where $s_2$ and $h_2$ are smooth functionals of the mean-field distribution. Next we present two important results on the contagion of divorce.
\begin{thm}  
In a short horizon, Breaking Up is Hard to Do, Unless Everyone Else is Doing it Too.
\end{thm}

Consider a society where the majority of marriages are stable. This means that both ${m}(t)$ and $h_2({m}(t))\geq 0 $ are high enough compared to $\underline{x}.$ Thanks to the contribution of $h_2$ in the feeling state dynamics, the generic dynamics given in (\ref{modelbook2}) has an higher value. The more the value of $h_2$ is, the more the feeling state will be. This means that
 $x(t,{m}(t))$  will go  to a higher state than the case without mean field.  Thus, Breaking Up is Hard to Do in a short horizon even if there is no effort from the couple.

\begin{rem}[Societal benefit!]
If $h_2({m}(t)) >> 0$ then the marriage remains  maintained over time even if the couple effort is minimal. This is because the mean-field contribution plays  the role of a positive effort. The network of the couple is having a big positive influence on their feeling state.
\end{rem}
\begin{thm}
If the mean-field has a tendency to the divorce states then, there is a contagious phenomenon for divorce, i.e., starting from $x_0 \in (0,2)$ the couple state will degrade due to the influence of the mean-field toward a negative feeling state.
In particular, Breaking Up is not seen as a negative thing in that society because  the majority is Doing it Too. Furthermore, stabilizing a marriage will require more effort, and hence it will be more expensive.
\end{thm}

Let $h_2({m}(t)) << 0.$  The mean-field is concentrated to negative values below the divorce threshold range. The majority  has  a tendency to the divorce. Thus, the feeling state of a generic couple goes towards to  negative values if no effort is made. A high effort is required to bring back and maintain the feeling state to a satisfactory one. The threshold effort $e(t)$ to balance is $\frac{h(x(t)) - h_2({m}(t))}{a} > 0$ which requires some time and cost. 
In this configuration, a divorce is not seen a negative thing by the society because the majority of the society is in a divorced state.
The negative term from the society will required more effort (and hence more cost) for the  stabilization of a marriage.

%
%
%
%
%
{\color{black}
\section*{Social welfare}
In this section we aim to maximize the social welfare. It consists to solve the following mean-field control problem:
\begin{eqnarray}
\label{mean_field2}
(\mathcal{P}_{3})
\left\{
\begin{array}{lll}
\displaystyle{ \sup_{e(.)} P_{\text{\tiny{social}}}(m_0, m, e)}, \\
\displaystyle{\mbox{subject to the feeling state dynamics }(\ref{modelbook2}) and }\ m(t)=\mathcal{L}(x(t)).
\end{array}
\right.
\end{eqnarray}
where $m(t)=\mathcal{L}(x(t))$ is the distribution of $x(t).$ Note that the social welfare problem $(\mathcal{P}_{3})$
 is different than $(\mathcal{P}_{2})$ because now a control action $e$ affects immediately the distribution $m$.

\begin{eqnarray}
\label{mf_payofftt2}
\displaystyle{  P_{\text{\tiny{social}}}(m_0, m, e)= \mathbb{E}\left(\hat{g}(x(T),m(T)){+}\int_0^T \left[\hat{s}(x(t), m(t)){-}c(e(t))\right] dt \right).}
\end{eqnarray}
Define the function $\hat{H}(x,m,p,q)=\hat{s}(x, m){-}c(e)+ f p+\sigma q.$
The stochastic maximum principle of mean-field type yields 
$$dp=-[ \hat{H}_{x}(x,m,p,q)+\mathbb{E}\hat{H}_{x,m}(x,m,p,q)]dt + q dB$$ where $e^*$ is the maximizer of $\hat{H}.$ Note that the impact of the society on the functions $\hat{s}, \hat{h}$ and $\hat{g}$ one gets the same optimality equations as above. However, the presence of the mean-field term $m$ changes drastically the behavior of the adjoint process $p.$
}
\section*{Discussions: Illustrative Model of Marriage} \label{threesteadystate}
In this section we present a mathematical model of marriage for a specific  functional $h$ in (\ref{modelrey}). In fact, we study the dynamics introduced in the "Mathematics of marriage" book \cite{refty0} in which we take into consideration a control term $a e$ and a noise term. The dynamical system in the case where
\begin{equation}
\label{h_expression} h(x) = rx - b \cdot \tanh(c x),
\end{equation} 
is written as
{\color{black}
\begin{eqnarray}
\label{model_example}
\displaystyle{dx(t) = (-r x(t) + b\cdot \tanh(c x(t)) + a e(t)) dt + \sigma d{B}(t),}
\end{eqnarray}
with $r,a,b,c>0$. The choice of the function $\tanh(x) = \frac{\exp(x) - \exp(-x)}{\exp(x) + \exp(-x)}$  in (\ref{h_expression}) is widely justified  in  the "Mathematics of marriage" \cite{refty0} as well as in evolutionary game theory \cite{maynardsmith} as a resulting from the imitative logit dynamics (see \cite{learning}). The model is widely supported by many psychologist and sociologist authors (see \cite{refty0}). Note that the function   $h$ defined at (\ref{h_expression}) does not belong to the class of functions studied in \cite{senti}.
 
}
\noindent  We choose $b=c=1$. Then, the sentimental dynamics with noise becomes
\begin{eqnarray} \label{withnoiset}
 dx(t) = [-r x(t) +  \tanh( x(t)) + a e(t)] dt + \sigma d{B}(t),
\end{eqnarray}
where the effort $e(t)\geq 0$ is a control variable. The parameter $r>0$ refers to the type of the society where the couple lives.
We also perform numerical examples to illustrate the effectiveness of the theoretical results. To do this, let consider  the mathematical model (\ref{modelrey}) with $h$ defined in (\ref{h_expression}). The cost and the well-being functionals are chosen as
\begin{eqnarray*} \label{c_and_s}
\displaystyle{c(e) = \frac{1}{2}e^2\;\;\text{and}\;\; s(x) = \bar{s} - 10 \exp(-x)\;\;\text{with}\;\; \bar{s} \geq 10.}
\end{eqnarray*}
It is worth noticing that the above choices of $c$ and $s$ ensure the existence of solutions for the problems (\ref{opt_open_loop_dyn}), (\ref{stoc_opt_open_loop_dyn}) (\ref{opt_feedback_effort}) and (\ref{mf_system}).
%
%

\begin{figure}[htb]
	\centering
		\includegraphics[width=0.5\textwidth]{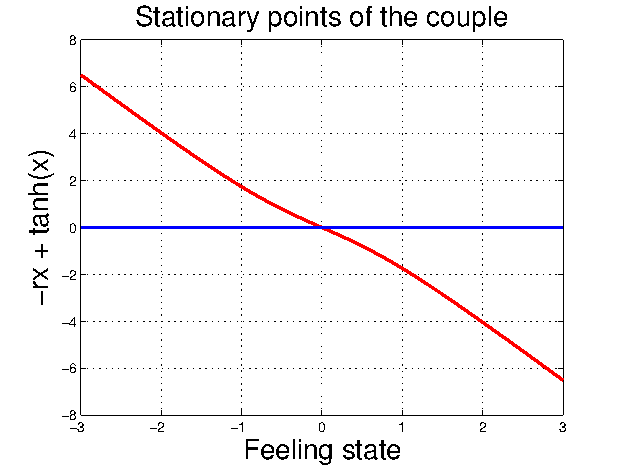}  
	    \caption{Steady states of the uncontrolled ($e=0$) dynamical system: high type society $r=2.5$ (one steady state).}
	\label{fig3:f1steadystate}
\end{figure}
\begin{figure}[htb]
	\centering
	    \includegraphics[width=0.5\textwidth]{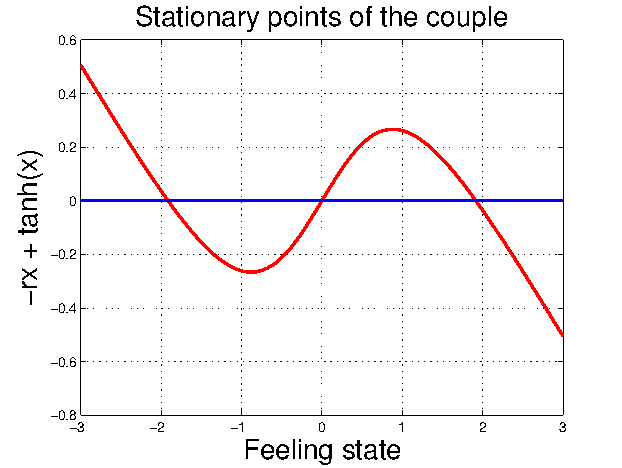}
	    \caption{Steady states of the uncontrolled ($e=0$) dynamical system: low type society $r=0.5$ (three steady states).}
	\label{fig4:f1steadystate}
\end{figure}
\subsection*{Steady states of the sentimental dynamics }
According to the type of the society, we analyze the uncontrolled ($e=0$) sentimental dynamics.
{\color{black}
\paragraph*{ In a high type society ($r\geq 1$).}
For $r\geq 1$ the uncontrolled system ($e = 0$) converges to zero independently of the starting point. This can be interpreted as if the couple feeling state will degrade over time if there is no effort. The uncontrolled feeling state has only one steady state which is illustrated in Figure \ref{fig3:f1steadystate}. 
\paragraph*{ In a low type society ($r<1$).} The degradation rate is less aggregative. We observe that in a society of low type, the uncontrolled system ($e = 0$) has  three steady states. Two of them are stable and one is unstable.
The unstable feeling state is the state $x = 0$ because the derivative at $0$ is $1-r>0$. The another two steady states are symmetric with respect to the origin. The one positive can be seen as if the marriage stays for life. A negative steady  state characterizes a divorce situation. The two non-zero states are stable. Indeed, the derivative at the steady points can be written as $-r+(1-\tanh^2(x^*) )= 1-r-(rx^*)^2$ which is negative when $rx$ is near $-1$ or $1.$ In figure \ref{fig4:f1steadystate} we illustrate the zeros of the drift function  $-rx +  \tanh(x)$ for different values of the societal parameter $r$.
%
%

}

\subsection*{ Evolution of uncontrolled feeling states }
Here, we show the societal influence on the sentimental dynamics of a  couple without effort ($e = 0$). In this case, we observe the variation in time of the feeling state for the two types of society. Numerical experiments for the deterministic ($\sigma = 0$) evolution and stochastic ($\sigma > 0$) dynamics in high type society ($r \geq 1$) and in low type society ($r < 1$) have been performed.
%
%

\begin{figure}[htb]
	\centering
		\includegraphics[width=0.5\textwidth]{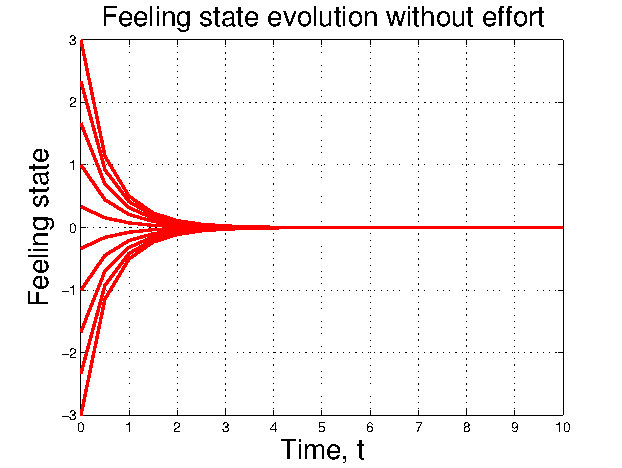}
	\caption{Uncontrolled sentimental dynamics without noise for $r=2.5$.}
	\label{fig5:f2sdet2nosigma}
\end{figure}

\begin{figure}[htb]
	\centering
		\includegraphics[width=0.5\textwidth]{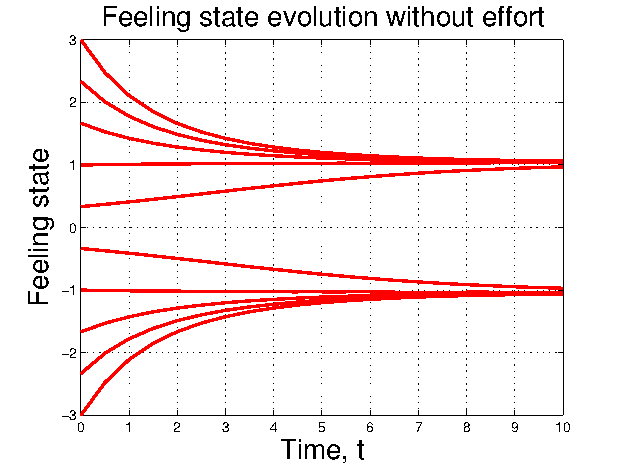}
	\caption{Uncontrolled sentimental dynamics without noise for  $r = 0.75$.}
	\label{fig6:f2sdet2nosigma}
\end{figure}

 A sample of 10 non-working couples with different initial feeling states are depicted in figure \ref{fig5:f2sdet2nosigma} and in figure \ref{fig6:f2sdet2nosigma}. We observe a fast convergence to the state $0$ in the high type case $(r=2.5).$ In the low type case $r=0.5$, we observe three steady states. Two (one positive, one negative) of them are stable states. 

\subsection*{ Stochastic evolution }
 Now we introduce a noise to better capture the fluctuations observed in real life. We observe that with  a small variance $\sigma$ the feeling state of the couple leads to a fluctuating trajectory around the deterministic sentimental dynamics. However, a bigger noise leads to two main branches and a non-zero probability to switch from one branch to another.

\begin{figure}[htb]
	\centering
		\includegraphics[width=0.6\textwidth]{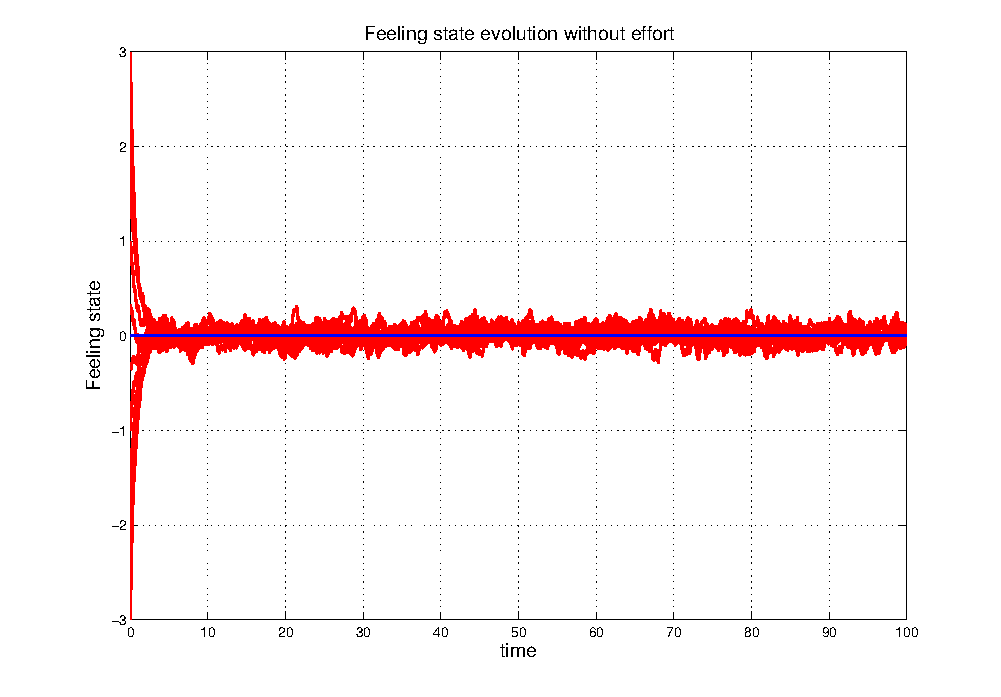}
	\caption{Noisy sentimental dynamics without control for $r=2$. }
	\label{fig7:f3sdet2withnoise}
\end{figure}

\begin{figure}[htb]
	\centering
		\includegraphics[width=0.6\textwidth]{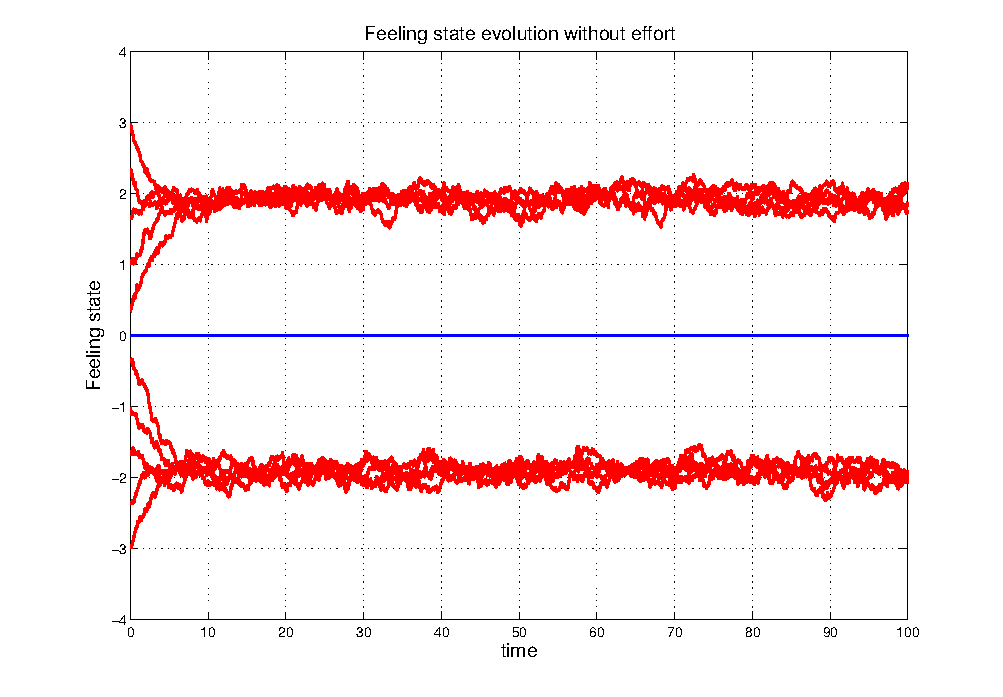}
	\caption{Noisy sentimental dynamics without control for $r= 0.5$. }
	\label{fig8:f3sdet2withnoise}
\end{figure}

Figure \ref{fig7:f3sdet2withnoise}   represents the stochastic evolution of  sentimental dynamics without control for $r=2$.  Figure \ref{fig8:f3sdet2withnoise} represents the noisy evolution of feeling states for low type $r= 0.5$.

%
%

\subsection*{Sentimental dynamics with effort }
\begin{figure}[htb]
	\centering
	\includegraphics[width=0.5\textwidth]{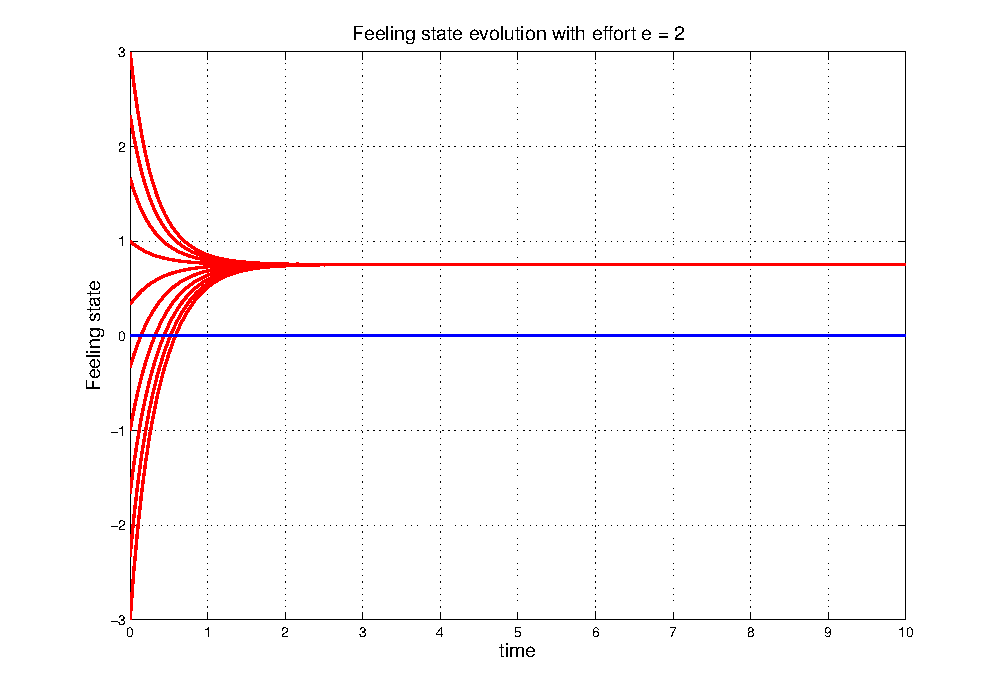}
\caption{Deterministic sentimental dynamics with control for  $r=3.5$.}
	\label{fig9:f4sdet2nosigmawithcontrol}
\end{figure}

\begin{figure}[htb]
	\centering
	\includegraphics[width=0.5\textwidth]{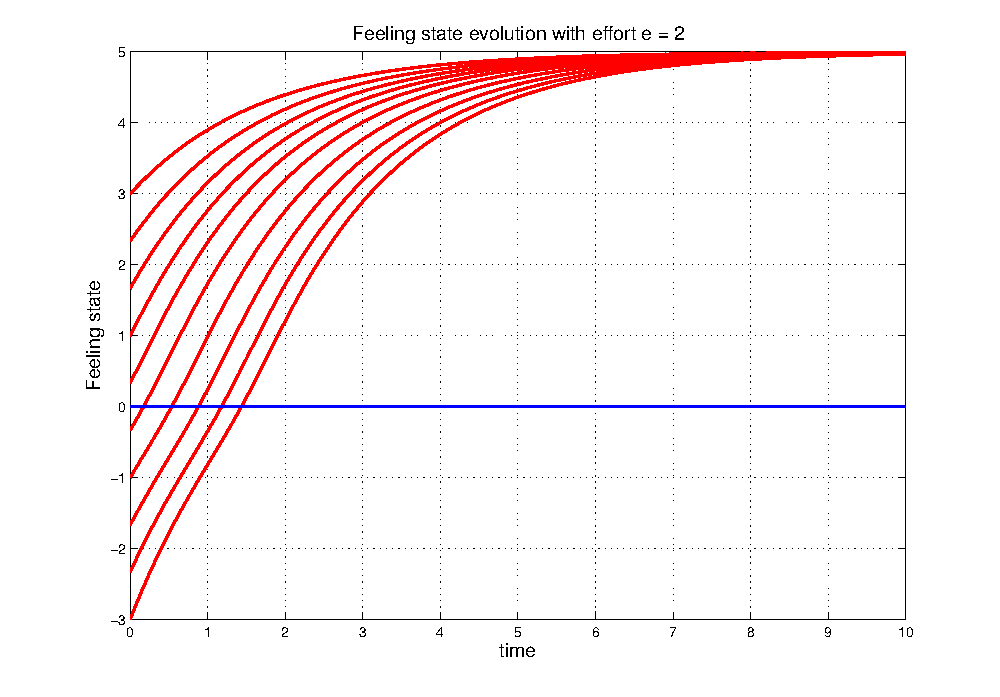}
\caption{Deterministic sentimental dynamics with control for  $r=0.6$.}
	\label{fig10:f4sdet2nosigmawithcontrol}
\end{figure}

\begin{figure}[htb]
	\centering
	\includegraphics[width=0.5\textwidth]{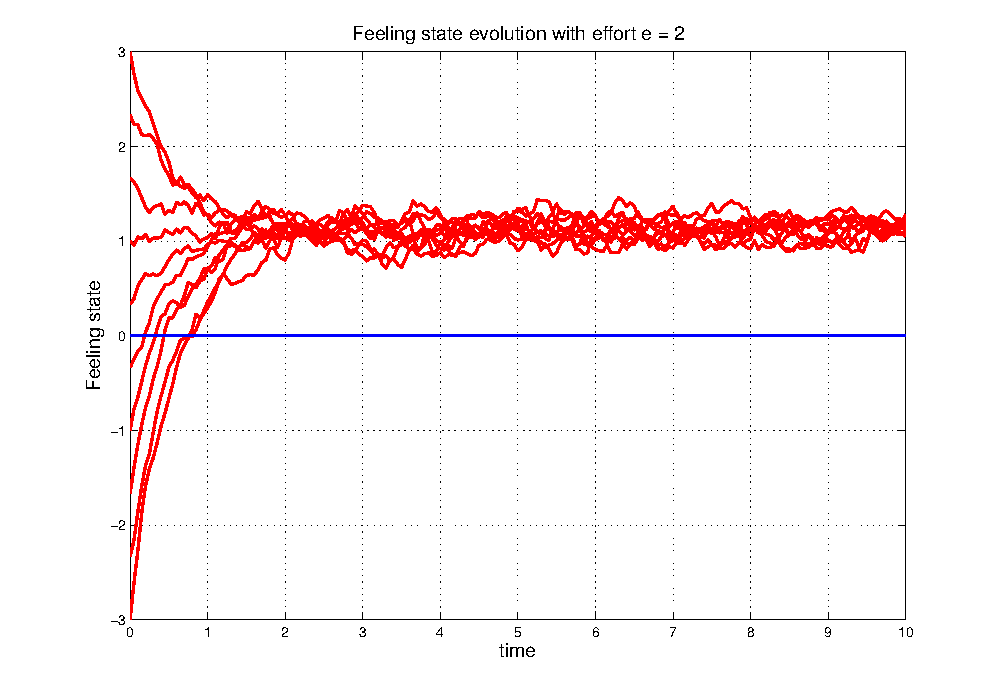}
	\caption{Stochastic feeling state with control  for $r=2.5$.}
	\label{fig11:f5sdet2withcontrol}
\end{figure}

\begin{figure}[htb]
	\centering
	\includegraphics[width=0.5\textwidth]{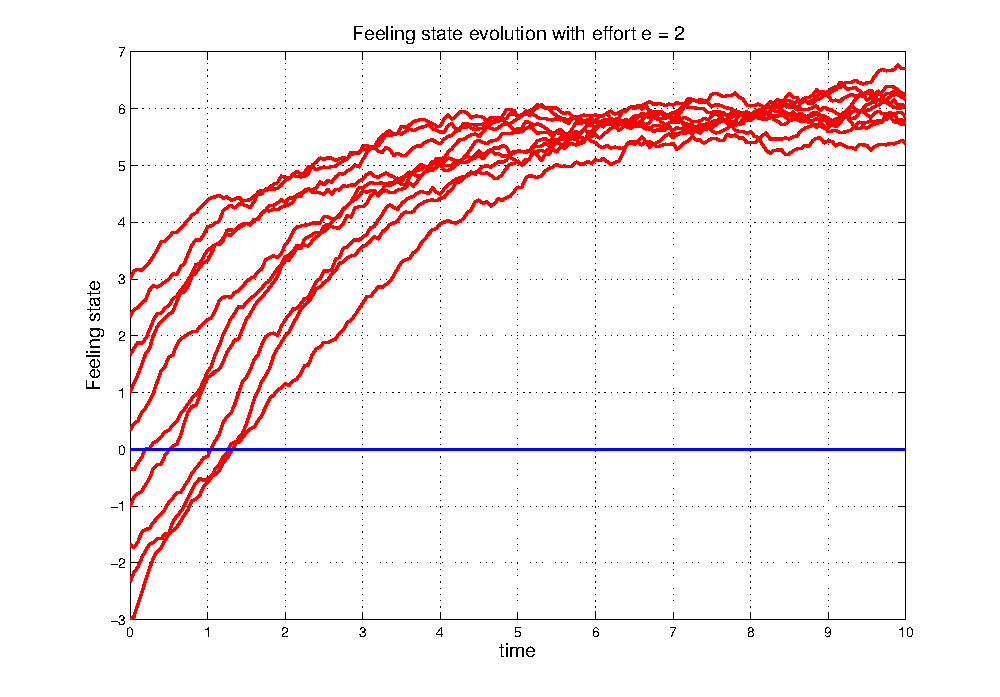}
	\caption{Stochastic feeling state with control  for $r=0.5$.}
	\label{fig12:f5sdet2withcontrol}
\end{figure}

Interestingly, when we introduce the effort, a catastrophic situation (starting from $x_0 <0$) can be reconstructed with effort and stabilize to a positive feeling state. These cases are illustrated in  Figures \ref{fig9:f4sdet2nosigmawithcontrol} and \ref{fig10:f4sdet2nosigmawithcontrol} for the deterministic dynamics and in Figures \ref{fig11:f5sdet2withcontrol}-\ref{fig12:f5sdet2withcontrol} for the stochastic dynamics. Thus, the effort plays an important role in maintaining the marriage or the cohabitation. On the other hand, high effort may be costly to the couple. Thus, it is crucial to find an acceptable tradeoff.

\begin{figure}[htb]
	\centering
\includegraphics[width=0.6\textwidth]{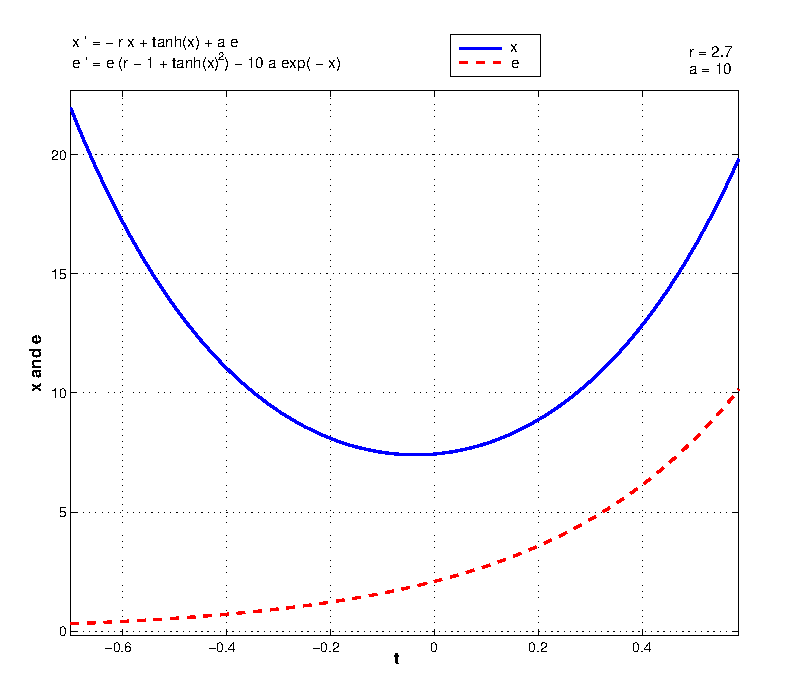}
	\caption{Open-loop optimal control system with high type $r=2.5$  as a function of time. }
	\label{fig13:f6optimalfeelinghigh}
\end{figure}

\begin{figure}[htb]
	\centering
\includegraphics[width=0.6\textwidth]{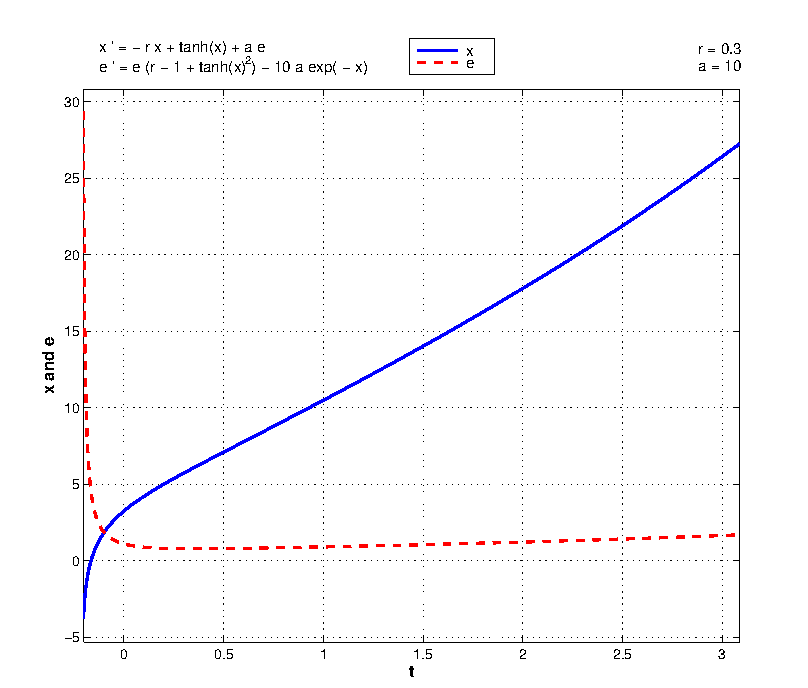}
	\caption{Open-loop optimal control system with low type $r=0.4$ as a function of time. }
	\label{fig14:f6optimalfeelinghigh}
\end{figure}

Next we depict the dynamics driven by optimal control as a function of time.
From the open optimal control of the non-noisy dynamics (\ref{opt_open_loop_dyn}), we observe two different situations. Figures \ref{fig13:f6optimalfeelinghigh} and \ref{fig14:f6optimalfeelinghigh} represent the trajectories of   $(x(t), e(t))$  for different initial conditions for both high and low type for $a=10$.

%
%

\begin{figure}[htb]
	\centering
   \includegraphics[width=0.5\textwidth,height=0.4\textwidth]{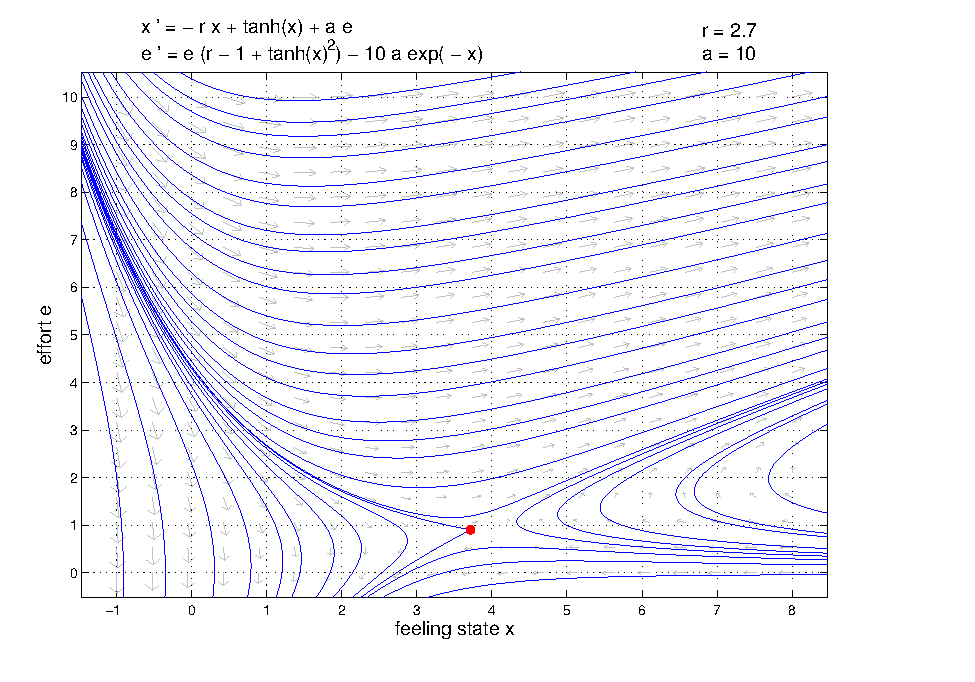}
	\caption{Vector field state versus effort for high type $r=2.5$.}
	\label{fig15:f7vectorfieldhigh}
\end{figure}

\begin{figure}[htb]
	\centering
   \includegraphics[width=0.5\textwidth,height=0.4\textwidth]{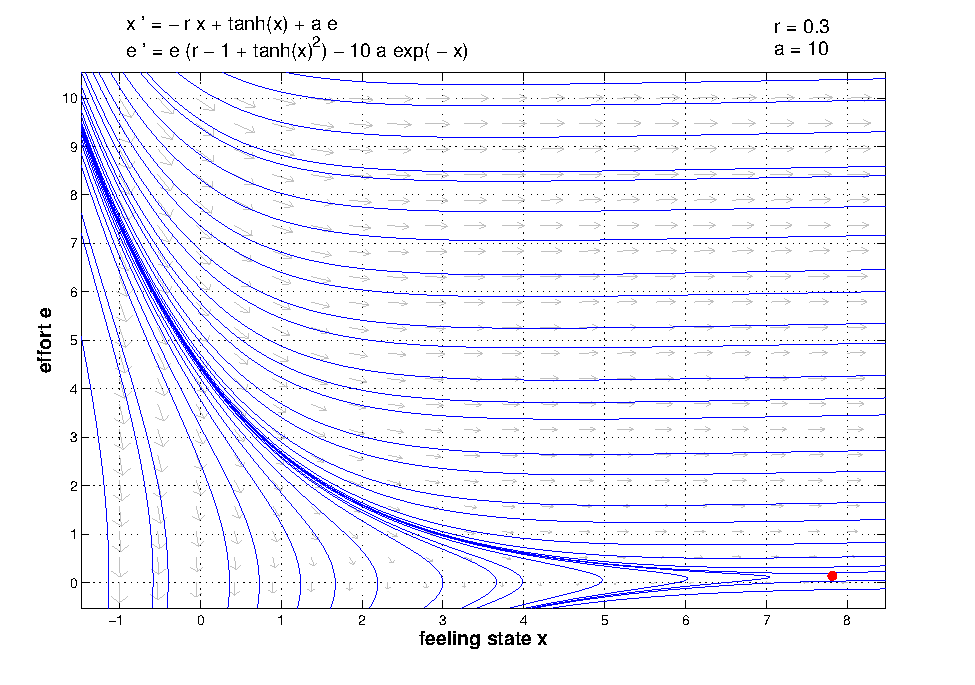}
	\caption{Vector field state versus effort for  low type $r=0.1.$}
	\label{fig16:f7vectorfieldhigh}
\end{figure}

Figure \ref{fig15:f7vectorfieldhigh}  and Figure \ref{fig16:f7vectorfieldhigh} represent the vector field between feeling state and effort of the couple for  high and low type. We observe there is an invariant set which can maintain very high effort and high feeling state.
%
%
%
%
%

\subsection*{ Stochastic optimal control }
\begin{figure}[htb]
	\centering
\includegraphics[width=0.7\textwidth,height=0.4\textwidth]{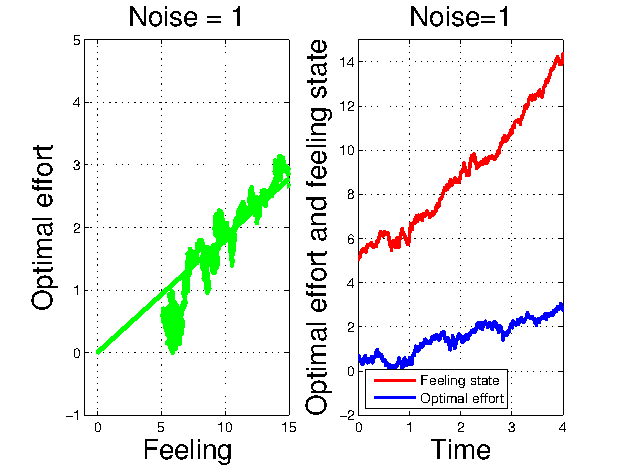}
	\caption{Open-loop stochastic optimal control system for a low type. }
	\label{fig17:noisyf6optimalfeeling}
\end{figure}

Figure \ref{fig17:noisyf6optimalfeeling} depicts the stochastic optimal control trajectories and optimal effort.  We observe that a small noise in the feeling state may have big consequences in the optimal effort (and the cost of effort).


%
%
%
%

\subsection*{Dynamic of  the optimal feedback effort }
The experiment consists in observing a couple over a finite horizon. We suppose that the couple finishes the horizon $T$ with the same heritage of feeling state whatever the initial feeling state is. In other terms, the functional $g$ is constant, therefore the terminal optimal effort is zero. To describe the dynamics of the optimal feedback strategy of such a couple, we applied an implicit (time) backward scheme  to the backward differential equation (\ref{opt_feedback_cost}). The Hamiltonian is approximated using Lax-Friedrichs approach to preserve monotonicity, consistency and continuity properties. Two numerical experiments are performed with respect to the type of the society.

\begin{figure}[htb]
\centering
\includegraphics[width=0.9\textwidth,height=0.5\textwidth]{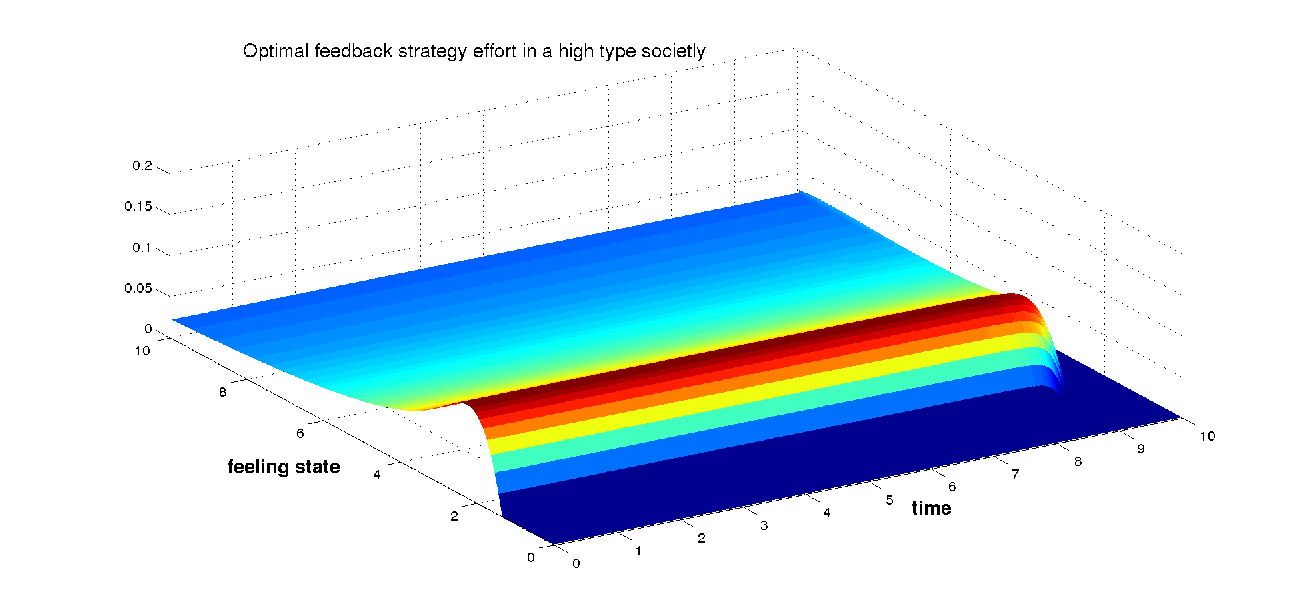}
\caption{The optimal strategy provided by a couple in a high type society $r=2.4$.}
\label{fig18_feedback_high}
\end{figure}

The couple starting with a low feeling state is  to provide no effort. For the conjoints of such couple, it is hard to resist due to the influence of the high type society. From Figure \ref{fig18_feedback_high}, one can  see that for a couple starting with a fairly good feeling state ($ \underline{x}+\epsilon \leq x \leq \underline{x}+K \epsilon $) the optimal strategy is to do actions. Indeed, the couple will break up because of the negative influence of society and will line up with a null optimal effort. Also, we observe the same phenomenon for a couple starting with a higher feeling state. The optimal effort is nevertheless important 
because starting at a higher feeling state  can counterbalance the negative influence of the society.

\begin{figure}[htb]
\centering
\includegraphics[width=0.9\textwidth,height=0.5\textwidth]{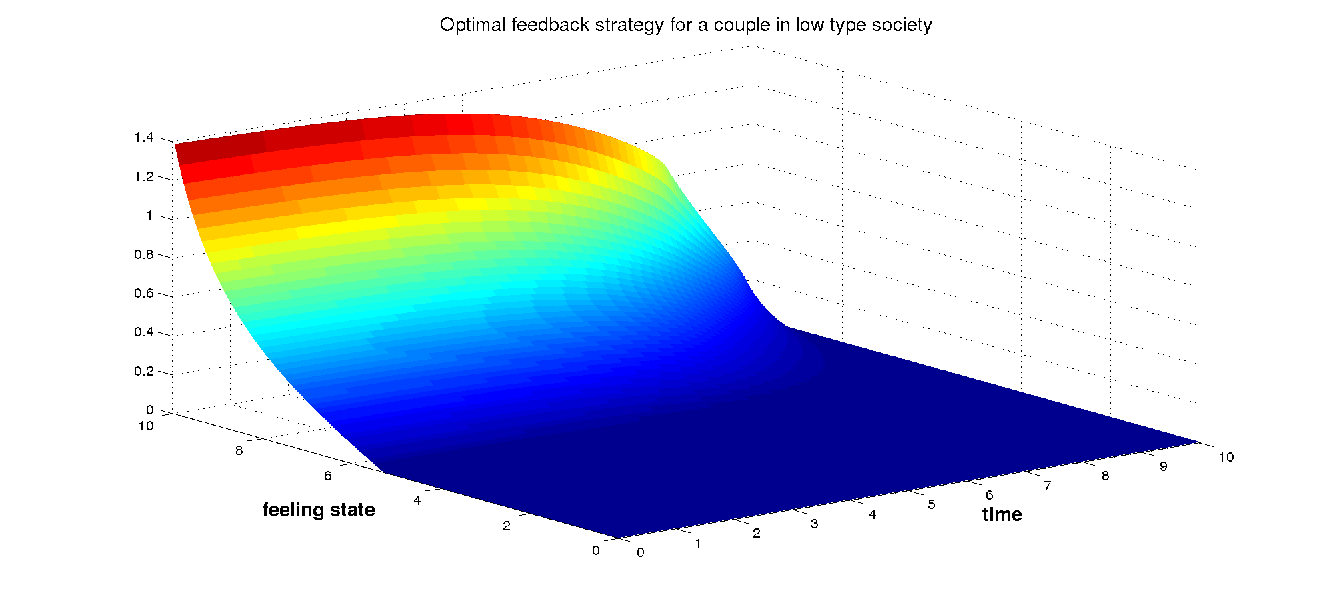}
\caption{ The optimal effort strategy of a couple in a low type society  $r=0.5$.}
\label{fig19_feedback_low}
\end{figure}

The optimal strategy for a couple starting with low feeling state is to  provide zero effort. This is not surprising since in a low type society living, a couple follows a societal phenomenon. We can see from Figure \ref{fig19_feedback_low} that the optimal strategy effort for a couple starting with high feeling state is a time-decreasing function. In that case, doing good action is motivated more by  seeking  ideal happiness than by  up-holding  the couple. Overall, we observe that the optimal strategy is cheaper in a low type society than in a high type society.
%
%
%

\subsection*{Mean field equilibrium trajectories}
We now address the numerical simulations for the mean-field equilibrium (\ref{mf_system}). The functional $\hat{h}$ and $\hat{s}$ are chosen as defined in (\ref{mf_h_and_s}), the functional $h$ is specified in (\ref{h_expression}), the effort cost and the well being  intrinsic to the couple are given by (\ref{c_and_s}).

The numerical experiment consists here in observing the evolution of the distribution $m$ representing the society in which a given couple lives. For the sake of simplicity, we suppose that the reference couple has an almost constant optimal effort in a particular window. This implies that the value  $v$ is linearized with respect to $x$ and constant in time. The numerical task is then to compute the Fokker-Planck-Kolmogorov equation in (\ref{mf_system}). The developed algorithm focuses on unnormalized distribution (positive measure).

\begin{figure}[htb]
\centering
\includegraphics[width=0.9\textwidth,height=0.5\textwidth]{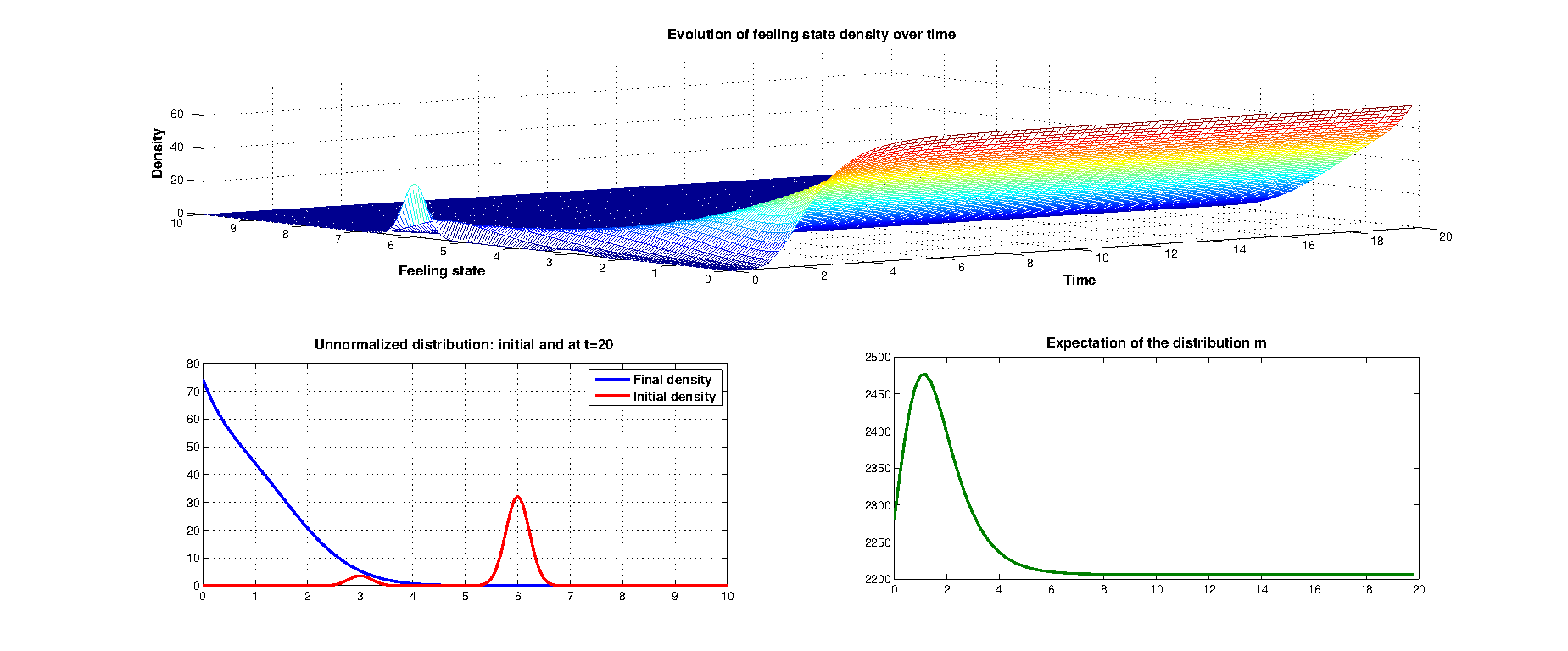}
\caption{Low type: Evolution of mean-field starting from Gaussian concentrated at $x=6$  with standard deviation $1.5$.}
\label{fig20mfe}
\end{figure}

\begin{figure}[htb]
\centering
\includegraphics[width=0.95\textwidth,height=0.65\textwidth]{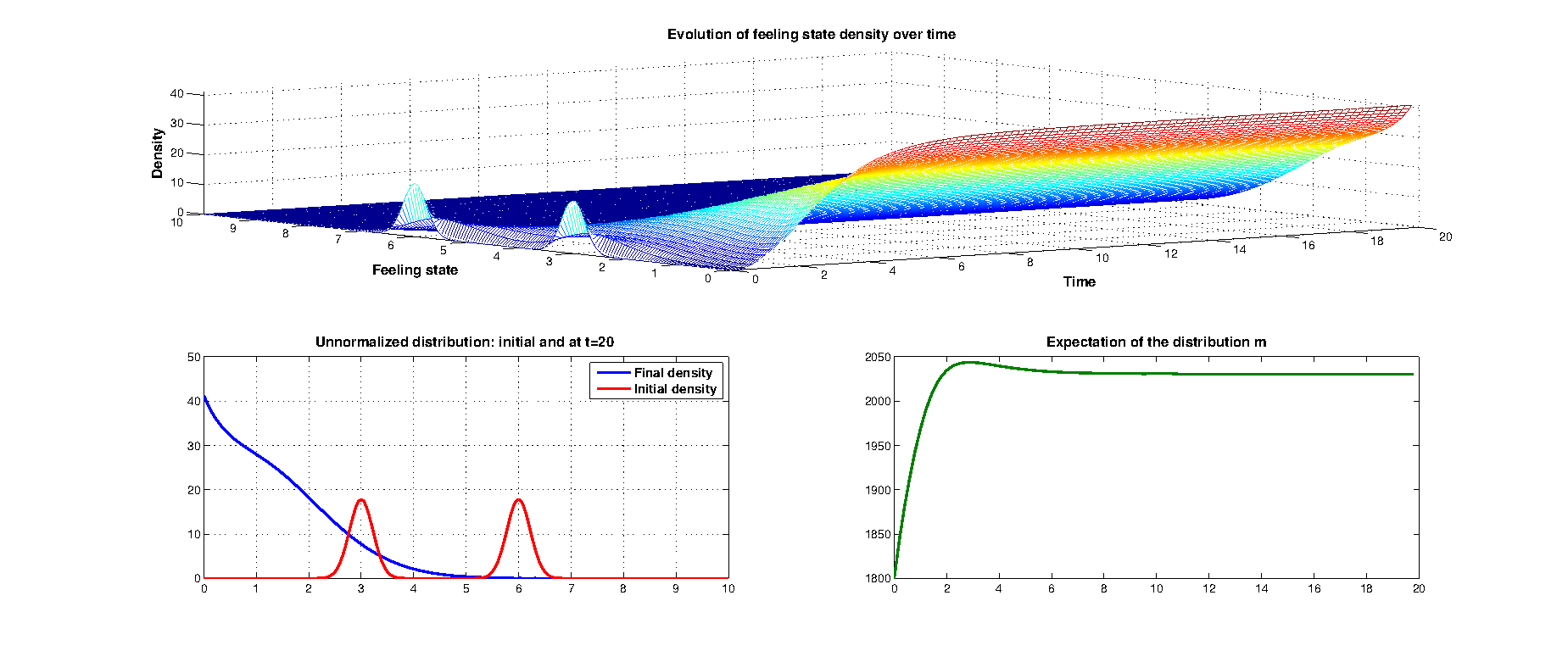}
\caption{Low type: Evolution of mean-field starting from mixture Gaussian concentrated at $x=3$ and $x=6$  with standard deviation $1.5$.}
\label{fig21mfe}
\end{figure}

Figure \ref{fig20mfe} represents the evolution of the mean-field equilibrium  with initial  Gaussian distribution centered at $x=6$  with standard deviation $1.5$ for a low type.
When the initial distribution is a mixture of Gaussian distributions  centered at $x=3$ and $x=6$  with standard deviation $1.5$, we observe in Figure \ref{fig21mfe} that  the state distribution will propagate rapidly to be concentrated at the extreme in a relatively short time.

\begin{figure}[htb]
	\centering
		\includegraphics[width=0.95\linewidth]{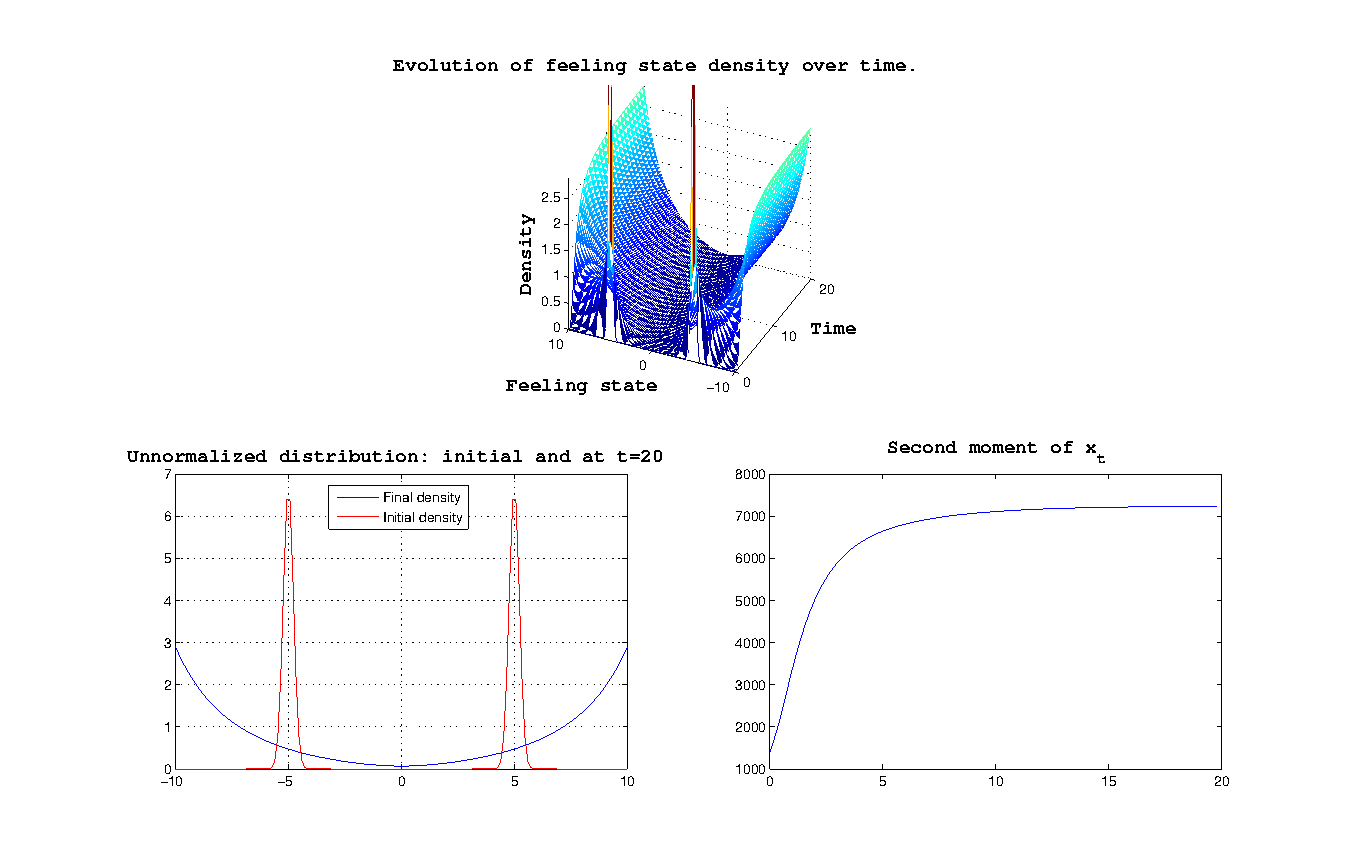}
	\caption{Mean-field equilibrium: the unnormalized distribution is concentrated at the two extreme boundaries.}
	\label{fig22mfet}
\end{figure}

For the mean-field equilibrium, the symmetric steady states observed previously in the case where the control was a constant, are not steady state anymore because the feedback control $e(t,x)$ is now dynamic in state as  time goes.
We observe that if  the initial mean-field is a mixture of Gaussians centered at $-5$ and  at $+5$  with standard deviation $0.05,$ then the feeling state will propagate and will be concentrated at the two extremes at the final time (Figure \ref{fig22mfet}).

\section*{Conclusion}
{\color{black} In this paper we have proposed a mean-field game model for marriage, cohabitation, divorce and remarriage. Our study suggests that the optimal effort may help in sustaining marriage and cohabitation if the cost of that effort is acceptable and the initial couple state is not too low. This helps us to understand at least theoretically the key processes related to marital dissolution, co-habitation separation, stability and fluctuation. The study can be useful to design or evaluate an adequate intervention.
It also suggests that knowing many divorced people may influence the status of a marriage. The more divorced people you know, the riskier your own marriage. However, a marriage doesn't break down just because friends are divorcing. Marital breakdown  depends on many factors including effort and mean field.
}
\section*{Acknowledgments}

We thank Alain Bensoussan, Galina Schwartz, Tatiana Sainati, Jose-Manuel Rey  and two anonymous reviewers for their constructive comments. 
%

%



\newpage

\end{document}